\documentclass{amsart}

\usepackage{amsthm, amssymb,amsrefs}

\newtheorem*{hyp*}{Hypothesis \((\ast )\)}
\newtheorem{thm}{Theorem}[section]

\newtheorem{lem}[thm]{Lemma}

\newtheorem{thmx}{Theorem}
\newtheorem{corx}[thmx]{Corollary}
\renewcommand{\thethmx}{\Alph{thmx}}

\newcommand{\irr}[1]{\text{Irr}(#1)}
\newcommand{\cod}[1]{\text{cod}(#1)}
\newcommand{\cd}[1]{\text{cd}(#1)}
\renewcommand{\c}[1]{{c}(#1)}
\renewcommand{\ker}[1]{\text{ker}(#1)}

\begin{document}

\title[Character codegrees of maximal class \(p\)-groups]{Character codegrees of maximal class $p$-groups} 

\author{Sarah Croome}
\address{%
Department of Mathematical Sciences\\
Kent State University\\
Kent, OH 44242\\
United States}

\email{scroome@kent.edu}

\author{Mark L. Lewis}
\address{Department of Mathematical Sciences\\
Kent State University\\
Kent, OH 44242\\
United States}
\email{lewis@math.kent.edu}


\begin{abstract} Let \(G\) be a \(p\)-group and let \(\chi\) be an irreducible character of \(G\). The codegree of \(\chi\) is given by \(|G:\text{ker}(\chi)|/\chi(1)\). If \(G\) is a maximal class \(p\)-group that is normally monomial or has at most three character degrees then the codegrees of \(G\) are consecutive powers of \(p\). If \(|G|=p^n\) and \(G\) has consecutive \(p\)-power codegrees up to \(p^{n-1}\) then the nilpotence class of \(G\) is at most 2 or \(G\) has maximal class.\end{abstract}

\maketitle

\section{Introduction}

In this paper all groups are finite. We investigate the set of codegrees of maximal class \(p\)-groups for prime \(p\).  For an irreducible character \(\chi\) of a group \(G\), the degree of \(\chi\) is given by \(\chi(1)\), and the kernel of \(\chi\) is denoted \(\ker{\chi}\). The codegree of \(\chi\) is defined as \(\cod{\chi}=|G:\ker{\chi}|/\chi(1)\), and the set of codegrees of the irreducible characters of \(G\) is denoted \(\cod{G}\). This definition for codegrees first appeared in \cite{OG}, where the authors use a graph-theoretic approach to compare the structure of a group with its set of codegrees. In more recent work \cite{codandnil}, Du and Lewis use codegrees to bound the nilpotence class of \(p\)-groups. Notably, they find that when a group has exactly 3 codegrees, its nilpotence class is at most 2. Their work motivates our investigation into the opposite extreme, when \(G\) is a maximal class \(p\)-group. Our first result demonstrates a surprising connection between these extremes: if \(G\) has order \(p^{n}\) and \(\cod{G}\) contains every power of \(p\) up to \(p^{n-1}\), then the nilpotence class of \(G\), \(c(G)\), can either be quite small or as large as possible, but not in between. This is the result of Theorem \ref{max_class_induction}.

\begin{thmx}
\label{max_class_induction}
If \(G\) is a \(p\)-group for some prime \(p\) such that \(|G|=p^{n}\ge p^2\), then \(\cod{G}=\{p^i\mid 0\le i\le n-1\}\) if and only if one of the following occurs:
\begin{enumerate}
\item \(G\cong \mathbb{Z}_{p^{n-1}}\times \mathbb{Z}_p\),
\item \(G \cong \mathbb{Z}_{p^{n-1}} \rtimes \mathbb{Z}_p\), \(G\) has nilpotence class 2, and \(n\ge 4\),
\item \(G\) has maximal class and \(|\cd{G}|=2\).
\end{enumerate}
\end{thmx}

 While examining maximal class \(p\)-groups, it became apparent that the codegrees of these groups are often consecutive powers of \(p\). From Theorem \ref{max_class_induction}, we know that a maximal class group with order \(p^{n}\) will only have all powers of \(p\) up to \(p^{n-1}\) as codegrees if the group has exactly two character degrees. When the group has three character degrees including \(p\), the codegrees are still consecutive powers of \(p\), but the largest power may vary.

\begin{thmx}
\label{max_class_3cd}
Let \(G\) be a maximal class \(p\)-group for some prime \(p\) such that \(\text{cd}(G)=\{1,p,p^b\}\). If \(|G|=p^n\), then \(\cod{G}=\{p^i\mid 0\le i\le c\}\), for some integer \(n-b\le c\le n-2\). 
\end{thmx}

Any metabelian maximal class group can have at most three character degrees, so from Theorem \ref{max_class_3cd} we know that the codegrees of such a group will be consecutive powers of \(p\). In this case, however, there are only two possibilities for the largest power of a codegree.

\begin{corx}
\label{max_class_metabelian_consecutive}
Let \(G\) be a metabelian maximal class \(p\)-group for some prime \(p\). If \(|G|=p^{n}\), then \(\cod{G}=\{p^i\mid 0\le i\le c\}\) where \(c=n-1\) or \(n-2\). 
\end{corx}

If \(G\) is a maximal class \(p\)-group that is also normally monomial, the powers of the codegrees will again be consecutive. This result is similar to that of Theorem \ref{max_class_3cd}, but does not depend on the number of character degrees of the group. 

\begin{thmx}
\label{normally_monomial}
Let \(G\) be a normally monomial maximal class \(p\)-group for some prime \(p\). Let \(|G|=p^{n}\), \(b(G)=\text{max}\{\chi(1)\mid \chi\in\irr{G}\}\), and \(b=\log_p(b(G))\). Then \(\cod{G}=\{p^i\mid 0\le i\le c\}\) for some integer \(c\ge n-b\).
\end{thmx}

It is easily shown that \(p^2\) and \(p^3\) are always included among the codegrees of maximal class \(p\)-groups (when the order of the group is large enough for the codegree to occur, e.g. \(p^3\in\cod{G}\) when \(|G|\ge p^4\)). If the order of \(G\) is at least \(p^6\), \(G\) also has \(p^4\) as a codegree.  This evidence, including Theorems \ref{max_class_induction}, \ref{max_class_3cd}, \ref{normally_monomial}, and Corollary \ref{max_class_metabelian_consecutive}, leads us to ask whether the codegrees of maximal class \(p\)-groups are always consecutive powers of \(p\).

\section{Main Results}

We begin with some necessary notation. The set of irreducible characters of a finite group \(G\) is written \(\irr{G}\), and the degrees of the irreducible characters of \(G\) form the set \(\cd{G}=\{\chi(1)\mid \chi\in\irr{G}\}\). In this paper, we consider only finite \(p\)-groups, and the degree of an irreducible character of such a group  is always a power of \(p\). Characters of degree \(1\) are called linear, with \(\text{Lin}(G)\) denoting the set of linear characters of \(G\). For a normal subgroup \(N\) of \(G\), there is a natural bijection between irreducible characters of  \(G/N\) and the irreducible characters of \(G\) whose kernels contain \(N\). We will therefore freely consider  \(\chi\in\irr{G/N}\) as a character of \(G\) whose kernel contains \(N\), and \(\ker{\chi}\) will usually refer to the kernel of \(\chi\) as a character of \(G\). Additionally, the codegree of \(\chi\) as an irreducible character of \(G/N\) is the same as the codegree of \(\chi\) as an irreducible character of \(G\), giving \(\cod{G/N}\subseteq \cod{G}\). The upper central series of \(G\) is the series \(1\le Z=Z_1\le Z_2\le \cdots \le Z_n\le G\), where \(Z(G/Z_i)=Z_{i+1}/Z_i\). The lower central series of \(G\) is \(1\le G_k\le \dots \le G'=G_2\le G_1=G\), where \(G_{i+1}=[G_i,G]\). For a \(p\)-group, the largest integer \(k\) such that \(G_k\) is nontrivial is called the nilpotence class of \(G\). 

The following lemma will be used repeatedly throughout this paper. It is a direct consequence of Proposition 2.5 in \cite{PPO1}, and It\^o's Theorem \cite[Theorem 6.15]{thebook}. 

\begin{lem}
\label{faithfulp}
If a finite \(p\)-group \(G\) has a faithful irreducible character of degree \(p\), then \(G\) has a normal abelian subgroup of index \(p\) and \(\cd{G}=\{1,p\}\).
\end{lem}

Notice the restriction \(n\ge 4\) in Theorem \ref{max_class_induction} (2). If \(G\) has class 2 and \(|G|=p^3\), then \(G\) need not be isomorphic to \(\mathbb{Z}_{p^2}\times \mathbb{Z}_p\), instead, \(G\) will be as in (3). In this case, note that \(G\) is extraspecial. 

\newpage

\begingroup
\def\thethmx{\ref{max_class_induction}}
\begin{thmx}
If \(G\) is a \(p\)-group such that \(|G|=p^{n}\ge p^2\) for some prime \(p\), then \(\cod{G}=\{p^i\mid 0\le i\le n-1\}\) if and only if one of the following occurs:
\begin{enumerate}
\item \(G\cong \mathbb{Z}_{p^{n-1}}\times \mathbb{Z}_p\),
\item \(G \cong \mathbb{Z}_{p^{n-1}} \rtimes \mathbb{Z}_p\), \(G\) has nilpotence class 2, and \(n\ge 4\),
\item \(G\) has maximal class and \(|\cd{G}|=2\).
\end{enumerate}
\end{thmx}
\endgroup

\begin{proof}
Assume first that \(G=\mathbb{Z}_{p^{n-1}}\times \mathbb{Z}_p\). Following the notation of Problem 2.7 of \cite{thebook}, as \(G\) is abelian, it is isomorphic to \(\widehat{G}\), the group of irreducible characters of \(G\). If \(\lambda \in \irr{G}\), the order of \(\lambda \) in \(\widehat{G}\) is given by \(|G:\ker{\lambda}|=\cod{\lambda}\). Hence the orders of characters in \(\widehat{G}\) correspond to the codegrees of \(G\), and we have \(\cod{G}=\{p^i\mid 0\le i\le n-1\}\).

Next, assume that \(G\cong \mathbb{Z}_{p^{n-1}} \rtimes \mathbb{Z}_p\) and \(G\) has nilpotence class 2. Let \(\langle x\rangle \cong \mathbb{Z}_{p^{n-1}}\) and choose \(\varphi \in \text{Aut}(\mathbb{Z}_{p^{n-1}}) \) as an automorphism of order \(p\) such that \(\varphi (x)=x^{1+p^{n-2}}\). Taking the semidirect product as \(\{ (a,b) \mid a \in \langle x \rangle, b \in \langle \varphi \rangle \}\), with multiplication given by \((a_1,b_1)(a_2,b_2)=(a_1(a_2^{b_1}),b_1b_2)\), we have \(G\cong \langle (x,1),(1,\varphi)\rangle \). Notice that \([(x,1),(1,\varphi)]=(x^{p^{n-2}},1)\), which has order \(p\), and hence \(G'=\langle (x^{p^{n-2}},1)\rangle \) has order \(p\).  As \((x^p,1)\in Z\), and \(|G:\langle (x^p,1)\rangle|=p^2\), we have \(Z=\langle (x^p,1)\rangle\). Now \(G/G'\cong \langle (x,1)\rangle \times \langle (1,\varphi )\rangle / \langle (x^{p^{n-2}},1)\rangle \cong \mathbb{Z}_{p^{n-2}}\times \mathbb{Z}_p\), which implies \(\cod{G/G'}=\{p^i\mid 0\le i\le n-2\}\). Since \(Z\) is cyclic, \(G'\) is the unique subgroup of \(G\) of order \(p\), and hence any nonlinear character of \(G\) must be faithful. Let \(\chi\in\irr{G}\) be nonlinear. Since \(|G:Z|=p^2\), by Corollary 2.30 of \cite{thebook}, \(\chi(1)=p\). Hence \(\cod{\chi}=p^{n-1}\), and we have \(\cod{G}=\{p^i\mid 0\le i\le n-1\}\). 

Finally, assume that \(G\) has maximal class and \(|\cd{G}|=2\). The quotient \(G/Z_{n-3}\) has order \(p^3\) and class 2. A nonlinear irreducible character of \(G/Z_{n-3}\) must have degree \(p\), as the center of \(G/Z_{n-3}\) has index \(p^2\). Hence \(\cd{G}=\{1,p\}\).  For \(i\le n-3\), the quotient \(G/Z_i\) is non-abelian with a cyclic center, and hence has a faithful nonlinear irreducible character with codegree \(p^{n-i-1}\). From the principal character and Corollary 2.3 of \cite{codandnil}, we always have \(1\) and \(p\) in \(\cod{G}\). Thus \(\cod{G}=\{p^i\mid 0\le i\le n-1\}\).

For the forward direction, assume first  that \(G\) is abelian with \(\cod{G}=\{p^i\mid 0\le i\le n-1\}\). Since \(p^{n}\notin \cod{G}\), \(G\) cannot be cyclic. As \(G\) is isomorphic to \(\widehat{G}\), there is an element of order \(p^{n-1}\) in \(G\). Since \(G\) is not cyclic, we must have \(G\cong \mathbb{Z}_{p^{n-1}}\times \mathbb{Z}_p\).

Now suppose \(G\) has class 2, \(n\ge 4\), and \(\cod{G}=\{p^i\mid 0\le i\le n-1\}\). Let \(\chi\in\irr{G}\) have codegree \(p^{n-1}\). If \(\chi\) is linear, then the kernel of \(\chi\) equals \(G'\) and has order \(p\). By Lemma 2.27 of \cite{thebook}, \(G/G'\) is cyclic, which implies \(G/Z\) is cyclic, contradicting that \(G\) has class 2. Hence \(\chi(1)=p\) and \(\chi\) is faithful. By Lemma 2.27 and Theorem 2.31 of \cite{thebook}, \(Z\) is cyclic and has index \(p^2\) in \(G\). Hence \(G\) has two noncentral generators, which implies \(|G'|=p\). Since \(Z\) is cyclic, \(G'\) is the unique normal subgroup of \(G\) of order \(p\), so an irreducible character of \(G\) is nonlinear if and only if it is faithful. Thus \(\cod{G/G'}=\{p^i\mid 0\le i\le n-2\}\). As an abelian group, \(G/G'\) is isomorphic to \(\widehat{G/G'}\). Hence \(\widehat{G/G'}\) (and therefore \(G/G'\)), has an element of order \(p^{n-2}\). Since \(G/G'\) is not cyclic, we must have \(G/G'\cong \mathbb{Z}_{p^{n-2}}\times \mathbb{Z}_p\).

Write \(G=\langle a,b\rangle\), where \(a^{p^{n-2}}\) and \(b^p\) are elements of \(G'\), but no smaller power of either element is in \(G'\). Suppose \(a^{p^{n-2}}=1\). Since \(|G:Z|=p^2\) and \(G/Z\) is not cyclic, it must be elementary abelian. Thus \(a^p\in Z\), and since \(Z\) has a unique subgroup of order \(p\), \(G'\) is contained in \(\langle a^p \rangle\). As \(a\langle a^p\rangle\) and \(b\langle a^p \rangle \) have order \(p\), and \(\langle a\langle a^p\rangle, b\langle a^p\rangle  \rangle = G/\langle a^p\rangle \), we have \(|G/\langle a^p\rangle |=p^2\). On the other hand, since \(p^{n-2}\) is the smallest power of \(a\) that equals \(1\), the order of \(\langle a^p\rangle \) is \(p^{n-3}\). Hence \(|G:\langle a^p\rangle |=p^3\), a contradiction. Thus \(a^{p^{n-2}}\neq 1\), and \(a\) has order \(p^{n-1}\). Now \(|G:\langle a\rangle|=p\). Suppose there is no element \(g\in G-\langle a \rangle\) which has order \(p\). Then  \(G\) has a unique subgroup of order \(p\), implying \(G\) is either cyclic or a generalized quaternion 2-group. Since generalized quaternion 2-groups have maximal class and \(|G|\ge p^4\), either case contradicts that \(G\) has class 2. Hence such an element \(g\) exists, which shows that \(G\cong\langle a \rangle \rtimes \langle g\rangle \cong \mathbb{Z}_{p^{n-1}} \rtimes \mathbb{Z}_p\). 


\emph{Claim}. \emph{If \(|G|=p^{n}\) and \(\cod{G}=\{p^i\mid 0\le i\le n-1\}\), then \(G\) has maximal class or \(c(G)\le 2\).}

\emph{Proof of claim.} Induct on \(|G|\). Clearly the result holds when \(|G|\le p^4\), so let \(|G|=p^{n}\) and \(\cod{G}=\{p^i\mid 0\le i\le n-1\}\). Assume \(\c{G}>2\). If \(\chi\in\irr{G}\) such that \(\cod{\chi}=p^{n-1}\), then \(\chi\) is faithful and has degree \(p\). Hence \(Z\) is cyclic, and \(G\) has a unique normal subgroup of order \(p\), say \(N\). The kernel of any non-faithful irreducible character must contain \(N\), and hence such a character will also be a character of \(G/N\). By Lemma \ref{faithfulp},  \(\cd{G}=\{1,p\}\). If \(\mu\) is a faithful irreducible character of \(G\) then \(\mu\) is non-linear and has degree \(p\), so \(\cod{\mu}=p^{n-1}\). This shows that an irreducible character of \(G\) has codegree \(p^{n-1}\) if and only if it is faithful. 
Thus \(\cod{G/N}=\{p^i\mid 0\le i\le n-2\}\), and by the inductive hypothesis \(\c{G/N}\le 2\) or \(\c{G/N}=n-2\). 

Suppose \(G/N\) has class \(n-2\). Then \(G\) has either class \(n-2\) or \(n-1\). If \(\c{G}=n-1\), we are done, so assume \(\c{G}=n-2\). Then \(|G_{n-2}|=p^2\), \(Z=G_{n-2}\), and \(|G'|=p^{n-2}\). By Lemma 1.1 of \cite{PPO1}, \(p^{n}=|G|=p|Z||G'|=p^{n+1}\), a contradiction.

Now suppose \(\c{G/N}\le 2\). If \(G/N\) is abelian, then \(\c{G}\le 2\) and we are done, so let \(\c{G/N}=2\), and assume \(\c{G}=3\). This implies \(N=G_3\). Let \(\gamma\in \irr{G/G_3}\) such that \(\cod{\gamma}=p^{n-2}\). If \(\gamma\) is linear then \(|\ker{\gamma}/G_3|=p\), so \(\ker{\gamma}=G'/G_3\), and by Lemma 2.27 (d) of \cite{thebook}, \(G/G'\) is cyclic, which is impossible. Hence \(\gamma(1)=p\) and \(\gamma\) must be a faithful character of \(G/G_3\). Put \(Z(G/G_3)=X/G_3\) and notice that \(X/G_3\) is cyclic. Since \(Z/G_3\) and \(G'/G_3\) are both contained in \(X/G_3\), and \(X/G_3\) has a unique subgroup of each possible index, we must have \(Z<G'\). Recalling Lemma 1.1 of \cite{PPO1}, we have \(p^{n}=p|Z||G'|\) and \(p^{n-1}=p|X/G_3||G'/G_3|\), hence \(|X:Z|=p\) and \(G'=X\). Put \(X=\langle a,G_3\rangle\), so \(Z=\langle a^p,G_3\rangle \). Since \(Z\) is cyclic, we either have \(a^p\in G_3\) or \(G_3\le \langle a^p \rangle\). If \(a^p\in G_3\), then \(Z=G_3\). Again applying Lemma 1.1, this shows that \(|G|=p^4\), so \(G\) has maximal class. Suppose instead \(G_3 \le \langle a^p \rangle\). The faithful character \(\gamma \in \irr{G/G_3}\) has center \(Z(\gamma)=X=G'\). Since \(G/Z(\gamma)\) is abelian, we have \(|G:X|=\gamma(1)^2=p^2\). This shows that \(G'=Z_2\). By Lemma 1.1 of \cite{PPO1}, \(p^{n}=|G|=p|Z||G'|=p\cdot p^{n-3}\cdot p^{n-2}=p^{2(n-2)}\), which shows that in this case as well we have \(|G|=p^4\). This proves the claim.

We may now assume \(G\) has maximal class and \(\cod{G}=\{p^i\mid 0\le i\le n-1\}\). Let \(\chi\in\irr{G}\) have codegree \(p^{n-1}\). If \(\chi\) is linear, then \(|\ker{\chi}|=p\), so \(\ker{\chi}=Z\). Since \(Z(\chi)=G\) and \(Z(\chi)/\ker{\chi}\) is cyclic, this implies \(G/Z\) is cyclic, which is impossible. Hence \(\chi\) is nonlinear, and must be faithful with degree \(p\). By Lemma \ref{faithfulp}, \(\cd{G}=\{1,p\}\). \end{proof}

The next lemma considers maximal class groups with  \(p^2\) as a character degree. The result is similar to that of Theorem \ref{max_class_induction} (3), with the largest codegree now reduced to \(p^{n-2}\). 

\begin{lem}
\label{max_class_cd3_consecutive}
Let \(G\) be a maximal class \(p\)-group with \(|G|=p^{n}\) and \(\cd{G}=\{1,p,p^2\}\). Then \(\cod{G}=\{p^i\mid 0\le i\le n-2\}\).
\end{lem}

\begin{proof}
A maximal class \(p\)-group with an irreducible character of degree \(p^2\) must have order at least \(p^5\), since \(|G:Z|\ge p^4 \) by Corollary 2.30 of \cite{thebook}. Let \(G\) be a maximal class \(p\)-group such that \(|G|=p^5\) and \(\cd{G}=\{1,p,p^2\}\). Suppose \(\chi\in\irr{G}\) has codegree \(p^4\), and notice that \(\chi\) must be a faithful character with degree \(p\). Lemma \ref{faithfulp} implies \(\cd{G}=\{1,p\}\), a contradiction, so \(\cod{\chi}\le p^3\) for all \(\chi\in\irr{G}\).
Since \(\c{G}=4\), Theorem 1.2 and Lemma 2.4 of \cite{codandnil} imply that \(|\cod{G}|\ge 4\), and hence \(\cod{G}=\{1,p,p^2,p^3\}\).

Now let \(|G|=p^{n}\) where \(n\ge 6\). Suppose \(\cod{\mu}=p^{n-1}\) for some \(\mu\in \irr{G}\). If \(\mu \) is linear, then \(|\ker{\mu}|=p\), which is impossible since the kernel of a linear character contains \(G'\), and \(G\) is maximal class. Hence \(\mu(1)=p\) and \(\mu\) is faithful, a contradiction by Lemma \ref{faithfulp}. Let \(\chi \in \irr{G}\) be faithful. Then \(\chi(1)=p^2\) and \(\cod{\chi}=p^{n-2}\) is the largest codegree of \(G\). If \(\cd{G/Z}=\{1,p\}\), then by Theorem \ref{max_class_induction}, \(\cod{G/Z}=\{1,p,p^2,\dots ,p^{n-2}\}\). If \(\cd{G/Z}=\{1,p,p^2\}\), then by the inductive assumption \(\cod{G/Z}=\{p^i\mid 0\le i\le n-3\}\). In either case, \(\cod{G}=\{p^i\mid 0\le i\le n-2\}\).
\end{proof}

This lemma can be extended to the case when the character degree \(p^2\) is replaced by an arbitrarily large power of \(p\). In this case, however, the largest codegree may vary. 

\begingroup
\def\thethmx{\ref{max_class_3cd}}
\begin{thmx}
Let \(G\) be a maximal class \(p\)-group such that such that \(\text{cd}(G)=\{1,p,p^b\}\). If \(|G|=p^n\), then \(\cod{G}=\{p^i\mid 0\le i\le c\}\), for some integer \(n-b\le c\le n-2\). 
\end{thmx}
\endgroup

\begin{proof}
First, notice that \(n\ge 2b+1\), as the square of the degree of every irreducible character of \(G\) must divide \(|G:Z|\). If \(\cod{\chi}=p^{n-1}\) for \(\chi\in \irr{G}\), then \(\chi\) is a faithful character of degree \(p\). By Lemma \ref{faithfulp}, this contradicts \(|G|=3\). Hence \(p^{n-2}\) is the largest possible codegree of \(G\). We proceed by induction on \(|G|\), with the base case \(|G|=p^5\) (and hence \(b=2\)) established by Lemma \ref{max_class_cd3_consecutive}. Now let \(|G|=p^{n}\), and let \(\chi\in \irr{G}\) be faithful. If \(\chi(1)=p\), then by Lemma \ref{faithfulp}, \(\cd{G}=\{1,p\}\), a contradiction. Hence \(\chi(1)=p^b\) and \(\cod{\chi}=p^{n-b}\). 

There are two possibilities for \(\cd{G/Z}\), either \(|\cd{G/Z}|=2\), or \(G/Z\) and \(G\) have the same three character degrees. In the first case, we can apply Theorem \ref{max_class_induction} to get \(\cod{G/Z}=\{p^i\mid 0\le i\le n-2\}\). Since \(b\ge 2\), we have \(\cod{\chi}=p^{n-b}\le p^{n-2}\), so \(\cod{G}=\cod{G/Z}\) and we are done. In the second case, we have \(\cod{G/Z}=\{p^i\mid 0\le i \le c'\}\), where \(c'\ge n-1-b\). If \(c'=n-1-b\), then \(\cod{G}=\{p^i\mid 0\le i \le n-b\}\). If \(c'>n-1-b\), then \(c'\ge n-b\), and \(\cod{G/Z}=\cod{G}\).
\end{proof}

Corollary \ref{max_class_metabelian_consecutive} is an easy corollary of Theorem \ref{max_class_induction} and Lemma \ref{max_class_cd3_consecutive}. 

\begingroup
\def\thethmx{\ref{max_class_metabelian_consecutive}}
\begin{corx}
Let \(G\) be a metabelian maximal class \(p\)-group. If \(|G|=p^{n}\), then \(\cod{G}=\{p^i\mid 0\le i\le c\}\) where \(c=n-1\) or \(n-2\). 
\end{corx}
\endgroup

\begin{proof}
If \(G\) is a metabelian maximal class \(p\)-group, then \(G'\) is abelian and has index \(p^2\) in \(G\). By It\^o's Theorem, if \(\chi\) is an irreducible character of \(G\), then \(\chi(1)|p^2\). Hence \(|\cd{G}|=2\) or \(\cd{G}=\{1,p,p^2\}\), and we can apply Theorem \ref{max_class_induction} and Lemma \ref{max_class_cd3_consecutive}, respectively.
\end{proof}

In \cite{codandnil}, Du and Lewis prove that \(1\) and \(p\) are always present in the set of codegrees of a finite \(p\)-group \(G\). When \(G\) has maximal class, we also have \(p^2\in\cod{G}\). 

\begin{lem}
\label{maxclass}
If \(G\) is a \(p\)-group that has maximal class, then \(p^2\in\cod{G}\).
\end{lem}

\begin{proof}
Let \(G\) have order \(p^{n+1}\) and class \(n\). The quotient \(G/Z_{n-2}\) is an extra-special group of order \(p^3\). If \(\chi\in\irr{G/Z_{n-2}}\) is nonlinear, then \(\chi\) must be faithful. Also notice that \(\chi(1)=p\), since \(|G/Z_{n-2}:Z(G/Z_{n-2})|=p^2\). Hence \(\cod{\chi}=p^2\), and since \(\cod{G/Z_{n-2}}\subseteq \cod{G}\), we have \(p^2\in \cod{G}\). 
\end{proof}

Following this trend, we find that whenever \(G\) has maximal class and \(|G|\ge p^4\), we have \(p^3\in\cod{G}\). 

\begin{lem}
\label{maxclass_p_cubed}
If \(G\) is a maximal class \(p\)-group such that \(|G|\ge p^4\), then \(p^3\in\cod{G}\).
\end{lem}

\begin{proof}
Let \(|G|=p^n\) where \(n\ge 4\). Since \(G\) has maximal class, \(Z(G/Z_{n-4})\) is cyclic and hence there exists a faithful character \(\chi\in\irr{G/Z_{n-4}}\). As \(G/Z_{n-4}\) is not abelian, \(\chi(1)>1\), but also \(\chi(1)^2\le |G:Z_{n-3}|=p^3\), so \(\chi(1)=p\). Thus \(\cod{\chi}=p^3\), so \(p^3\in\cod{G/Z_{n-4}}\subseteq \cod{G}\).
\end{proof}

We would like to be able to continue increasing the order of the group to obtain each next largest codegree, however our results so far are limited to codegrees up to size \(p^4\), as in the following lemma.

\begin{lem}
\label{maxclass_p_4}
Let \(G\) be a maximal class \(p\)-group such that \(|G|\ge p^6\). Then \(p^4\in \cod{G}\). 
\end{lem}

\begin{proof}
Let \(|G|=p^n\) where \(n\ge 6\). Since \(|G:Z_{n-5}|=p^5\), the largest possible degree of an irreducible character of \(G/Z_{n-6}\) is \(p^2\). By Theorem \ref{max_class_induction}  and Lemma \ref{max_class_cd3_consecutive}, we must have \(p^4\in\cod{G/Z_{n-6}}\subseteq \cod{G}\). \end{proof}

This bound is sharp, and the groups listed in the the small group library of Magma \cite{magma} as SmallGroup\((3^5,i)\) for \(i=28,29,\) and \(30\) are examples of maximal class groups of order \(p^5\) with \(p^4\) not included among their codegrees. When \(|G|=p^{n}\) and \(c(G)=n-1\), it is not surprising that \(p^{n-1}\) is not always included in the set of codegrees, as this would imply \(G\) has a faithful character of degree \(p\), which in turn implies that \(\cd{G}=\{1,p\}\), and the existence of groups of maximal class with other sets of character degrees is well known. Whether this pattern continues is an interesting question. We have not been able to prove that \(p^5\) must be a codegree for maximal class groups of order at least \(p^8\), nor have we found a counterexample. If a counterexample \(G\) of order \(p^8\) exists, then \(G\) must have \(p^3\in\cd{G}\), but no faithful irreducible character of \(G\) can have degree \(p^3\).  

An irreducible character \(\chi\) of a group \(G\) which is induced from a linear character of a subgroup of \(G\) is called monomial. If this subgroup is normal in \(G\) then \(\chi\) is called normally monomial. If every irreducible character of \(G\) is monomial, we say \(G\) is an \(M\)-group, and if every \(\chi\in\irr{G}\) is normally monomial, then \(G\) is also said to be normally monomial. All \(p\)-groups are nilpotent, and by Corollary 6.14 of \cite{thebook}, they are therefore \(M\)-groups. Our final theorem is restricted to normally monomial maximal class \(p\)-groups. In a group of this type, a faithful irreducible character will always have the largest degree in \(\cd{G}\), allowing us to apply induction to \(G/Z\).

\begingroup
\def\thethmx{\ref{normally_monomial}}
\begin{thmx}
Let \(G\) be a normally monomial maximal class \(p\)-group. Let \(|G|=p^{n}\), \(b(G)=\text{max}\{\chi(1)\mid \chi\in\irr{G}\}\), and \(b=\log_p(b(G))\). Then \(\cod{G}=\{p^i\mid 0\le i\le c\}\) for some integer \(c\ge n-b\).
\end{thmx}
\endgroup

\begin{proof}
Let \(G\) be nonabelian with order \(p^3\) and let \(\chi\in\irr{G}\). Since \(\chi(1)^2\le |G:Z|=p^2\), we have \(\cd{G}=\{1,p\}\). If \(\chi(1)=p\), then by Lemma 2.1 of \cite{codandnil} \(\cod{\chi}\ge p^2 \). An irreducible character of \(G\) with codegree \(p^3\) must be faithful and linear, which is impossible since \(G\) is nonabelian. Hence \(\cod{\chi}=p^2\), and \(\cod{G}=\{1,p,p^2\}\).

Now let \(G\) be normally monomial with \(|G|=p^{n}\) and \(\c{G}=n-1\). Since quotients of normally monomial groups are also normally monomial, we can apply induction and assume \(\cod{G/Z}=\{p^i\mid 0\le i \le c'\}\), where \(c'\ge n-1-b'\) and \(b'=\log_p (b(G/Z))\). Notice that \(b'\le b\). If \(b'<b\), then \(b'\le b-1\), so \(c'\ge n-1-b' \ge n-1-(b-1)=n-b\). Hence \(\cod{G/Z}=\{p^i\mid 0\le i \le c'\}\subseteq \cod{G}\).

Let \(\chi\in\irr{G}\). If \(\chi\) is not faithful, then \(\chi\) can be equated with an irreducible character of \(G/Z\) and hence \(\cod{\chi}\in \cod{G/Z}\). Thus assume \(\chi\) is faithful. As \(G\) is normally monomial, there exists a subgroup \(H\unlhd G \) and \(\lambda \in \text{Lin}(H)\) such that \(\chi=\lambda^G\) and \(\chi(1)=\lambda^G(1)=|G:H|\). Since \(H'\) is characteristic in \(H\) which is normal in \(G\), we have \(H'\unlhd G\). Thus \(H'\le \text{core}_G(\ker{\lambda}) =\cap_{x\in G}(\ker{\lambda})^x\). By Lemma 5.11 of \cite{thebook}, this is equal to \(\ker{\lambda^G}=\ker{\chi}=1\), which shows that \(H'\) is trivial and hence \(H\) is abelian. If \(\theta\in \irr{G}\), then by It\^o's Theorem, \(\theta\) divides \(|G:H|=\chi(1)\). Thus \(\chi(1)=b(G)\), and \(\cod{\chi}=p^{n-b}\in\cod{G/Z}\).

Assume now that \(b'=b\), so \(\cod{G/Z}=\{p^i\mid 0\le i \le c'\}\), where \(c'\ge n-1-b\).  Let \(\chi\in \irr{G}\) and note again that if \(\chi\) is not faithful then \(\cod{\chi}\in\cod{G/Z}\) so we may assume \(\chi\) is faithful. Since \(G\) is normally monomial, we have (as in the previous paragraph) that \(b(G)=\chi(1)\), and therefore \(\cod{\chi}=p^{n-b}\). Thus either \(\cod{\chi}\in\cod{G/Z}\), in which case \(\cod{G}=\cod{G/Z}=\{p^i \mid 0\le i \le c\}\), where \(c\ge n-b\), or \(\cod{G}=\{1,p,\dots,p^{n-1-b},p^{n-b}\}\). \end{proof}


\end{document}